# The Polygons of Albrecht Dürer -1525
G.H. Hughes

The early Renaissance artist Albrecht Dürer published a book on geometry a few years before he died. This was intended to be a guide for young craftsmen and artists giving them both practical and mathematical tools for their trade. In the second part of that book, Durer gives compass and straight edge constructions for the 'regular' polygons from the triangle to the 16-gon. We will examine each of these constructions using the original 1525 text and diagrams along with a translation. Then we will use Mathematica to carry out the constructions which are only approximate, in order to compare them with the regular case. In Appendix A, we discuss Dürer's approximate trisection method for angles, which is surprisingly accurate. Appendix B outlines what is known currently about compass and straightedge constructions and includes two elegant 19th century constructions.

There is list of on-line resources at the end of this article which includes two related papers at DynamicsOfPolygons.org. The first of these papers is an analysis of the dynamics of the non-regular 'Dürer polygons' under the outer billiards map (Tangent map). The second paper addresses the general issue of compass and straightedge constructions using Gauss's *Disquisitiones Arithmeticae* of 1801.

# Part I - Dürer's life

Albrecht Dürer was born in Nuremberg (Nürnberg) in 1471. Nuremberg lies on the Pegnitz River and is surrounded by the farmlands and woods of Bavaria. Medieval Germany was part of the Holy Roman Empire and when Frederick III died in 1493, his son Maximillian I became ruler of the Roman Empire and Germany. Nuremberg was an important city in this empire because it lay on the trade routes between southern and northern Europe. Under Maximillian I, Nuremberg was regarded as the 'unofficial capitol' of the Roman Empire.

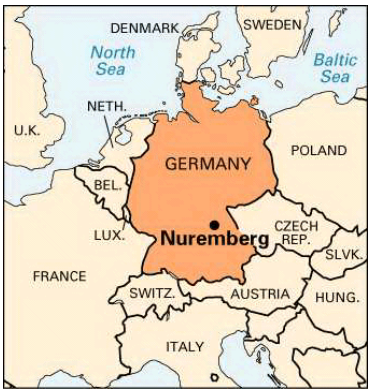

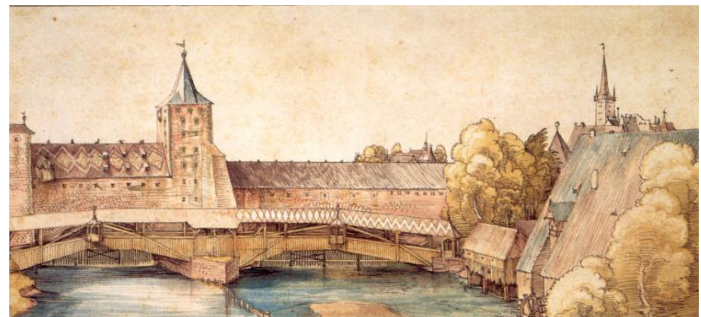

*The walkway street market in Nuremberg* by Albrecht Dürer

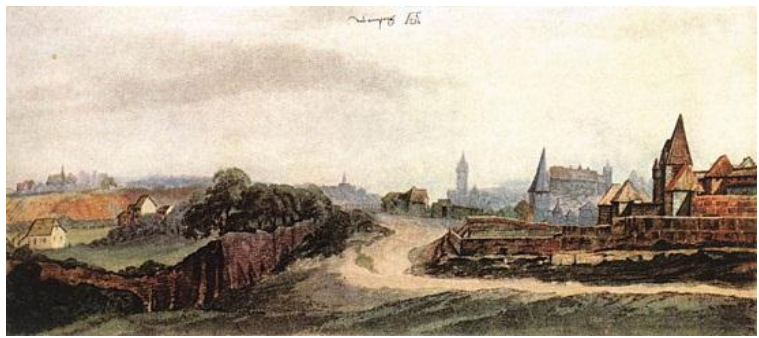

*View of Nuremberg from the West* by Albrecht Dürer

Nuremberg was famous for craftsmen who produced weapons, nautical instruments, jewelry, woodcuts and engravings. In 1450, in nearby Mainz, Johannes Gutenberg introduced a printing press using movable type and book printing soon became a thriving industry in Nuremberg. The printing press and the town of Nuremberg played a major part in the Reformation movement in Germany.

Dürer's father came to Nuremberg from his native Hungary in1455. He worked as a goldsmith in the shop of Hieronymus Holper and 12 years later married Holper's daughter. Albrecht was the third child of this marriage and within 20 years he had a total of 17 siblings – only three of which survived childhood. (Nuremberg suffered a number of plagues in the 15$^{th}$ and 16$^{th}$ centuries.)

According to Albrecht's book of family matters, his father had been trained in the Netherlands "*with the great Masters*". The Elder Dürer was industrious but not particularly prosperous. Following tradition, Albrecht became an apprentice to his father and learned the skills of a goldsmith. In medieval Germany, goldsmiths were closely aligned with engravers because the skills were very similar. In the words of Edwin Panofsky " *the greatest engravers of the 15$^{th}$ century were originally not painters or book illuminators, but goldsmiths…*"

It soon became clear that Albrecht had natural talent for drawing. Below is a self-portrait, done in silverpoint at age 13. A year later he became an apprentice to one of the most famous artists and printmakers in Nuremberg, Michael Wolgemut.

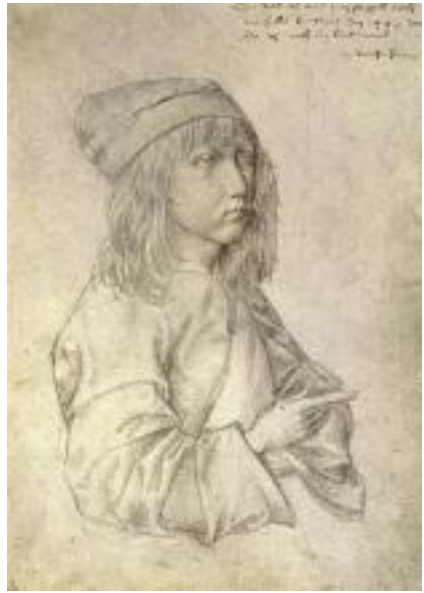

Upon completing his apprenticeship at age 19, Dürer following the tradition of young artists and traveled throughout Germany to study the work of master artists and printmakers. He traveled to Colmar and later Basel which was a center for bookmaking. Dürer's godfather was a prosperous publisher and the young Dürer already had influential connections in Basel. He did a woodcut portrait of St Jerome which appeared in 1492 in the Letters of St. Jerome published by Nicolaus Kessler. Other publishers in Basel admired Dürer's work and asked him to provide woodcuts for their books. From Basel he traveled to Strassburg but he was called home by his father in 1494. The Elder Dürer had arranged a marriage with Agnes Fry and the dowry would soon enable him to open his own shop.

After getting married in the summer of 1494, Dürer set off to Italy where his life-long friend Willibald Pirckheimer was studying in Padua. The warmth of Italy was a respite from the plague that had broken out in Nuremberg. Dürer witnessed the art of the Italian Renaissance and marveled at the level of respect showed toward artists. In his native Germany there were no 'artists' – only craftsmen. Dürer was determined that Germany would participate in the regrowth of science and the arts that he was witnessing in Italy.

The Renaissance was accompanied by newfound respect for geometry, mathematics and the sciences. In Venice, Dürer found artists such as Luca Pacioli who would publish *The Divine Proportion* in 1509 with drawings by Leonardo da Vinci (1452- 1519).  Dürer was not highly educated but he was familiar with the works of Euclid (c.300 BC) and Ptolemy (AD 90 – AD 168) and on a later trip to Venice,  Dürer bought a copy of Euclid for his library. Throughout his life, Dürer exchanged views on art and science with Pirckheimer and other Nuremberg scholars such as the mathematician Johannes Werner.

When Dürer returned to Nuremberg in 1495 he opened his own shop and soon became one of the most respected printmakers in Germany. He was an accomplished painter but his true passion was printmaking. The prints could be sold at a price which was affordable by the common person and Dürer made every effort to produce prints which were interesting and enlightening. Below is a watercolor of the courtyard of the former castle in Innsbruck, and on the right is a woodcut of three peasants in conversation.

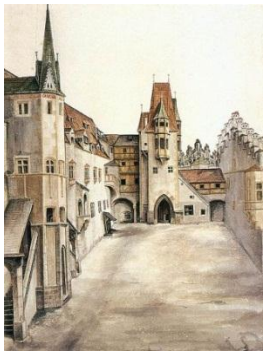 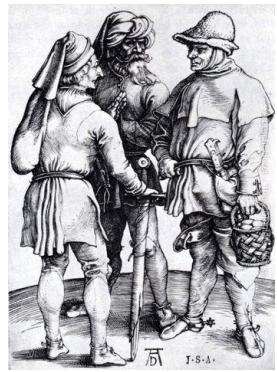

In 1505 Dürer made another trip to Italy. In general, the Italians regarded the Germans as poorly educated brutes, but they had a great respect for Dürer as an artist and a man. He became a close friend of Giovanni Bellini whom Dürer described as "*the best painter in Venice*". Dürer was becoming one of the most influential artists in Northern Europe. When he returned to Germany, the Emperor Maximillian I, honored him with commissions and later gave him a pension. In 1519 Dürer painted  a portrait of Maximillian holding a pomegranate which was his personal symbol.

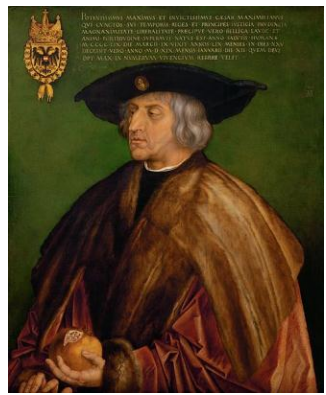

When Maximillian was succeeded by his grandson Charles V, Dürer traveled to the Netherlands to meet with Charles. He was honored throughout his trip, but an excursion to Zeeland in 1520 to view a beached whale, apparently led to a bout with malaria which contributed to his demise in 1528 at age 57. Dürer left behind a considerable estate and a legacy of more than a thousand drawings, hundreds of engravings and dry points, and about seventy paintings. He also wrote and illustrated books on geometry, city fortification, and human proportion.

In keeping with the tradition of the early Renaissance period , most of Dürer's work was religious in nature, but he also produced some notable secular works such as the famous *Melancholia I* copperplate engraving from 1514 (see below). This is a study in perspective and geometry which also contains a number of mathematical references such as the 'magic' square in the upper right corner. This square sums to 34 in rows, columns and diagonals. Note that the bottom row contains 15 and 14 which is the date of the engraving, and the numbers 4 and 1 may represent D and A - Dürer's monogram. In addition, pairs of numbers symmetrical to the center sum to 17. This magic square may have originated with Cornelius Agrippa's *De Occulta Philosophia*.

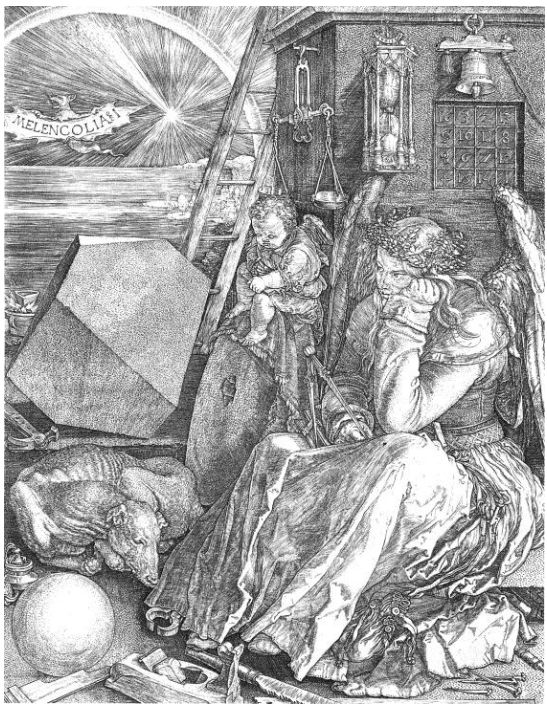
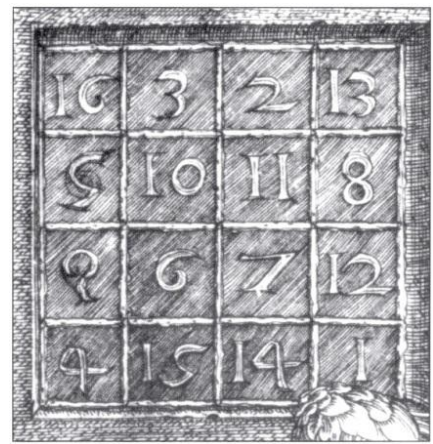

There is a strange polyhedron which appears to be a truncated rhomboid. It may have started out as a cube which was stretched and then truncated. Mathematica has a data-base of polyhedra which includes DuerersSolid (this is the accepted spelling when the ü is not desirable) - type, **PolyhedronData["DuerersSolid"].** The skeleton of Dürer's solid is a generalized Peterson graph. There is still doubt as to the actual geometry of the polyhedron and its meaning. Dürer developed innovative methods for constructing polyhedra using 'nets' and it is possible that the net for this polyhedron is related to the magic square.

In the lower left part of the engraving there are tools such as compass, scale and hourglass. These tools lay unused on the ground and there is endless speculation about the meaning of this. Dürer appears to have suffered from early bouts of arthritis and after visiting the Netherlands in 1521, he may have also suffered from malaria. In addition, his arranged marriage had been a disappointment and he had no children.

Late in his life, Dürer embarked on an ambitious project to convey to other artists and craftsmen the theoretical skills of their craft. Dürer felt that his fellow German artists were equal to all others in practical skill and power of imagination, but they were inferior to the Italians in rational knowledge. In keeping with the spirit of the Renaissance Dürer felt that "*geometry is the right foundation of all painting*".  But Dürer was also a practical man and he later wrote: "*As for geometry it may prove the truth of some things; but with respect to others we must resign ourselves to the opinion and judgment of men.*"

Dürer finished just the first stages of his project: The Four Books of Measurement with Compass and Ruler (*Underweysung der Messung mit dem Zirckel und Richtscheyt)* which was published in 1525 and the Four Books of Human Proportion (*Vier Buchen von Manschlicher Proportion*) which was published in 1528 just after his death.

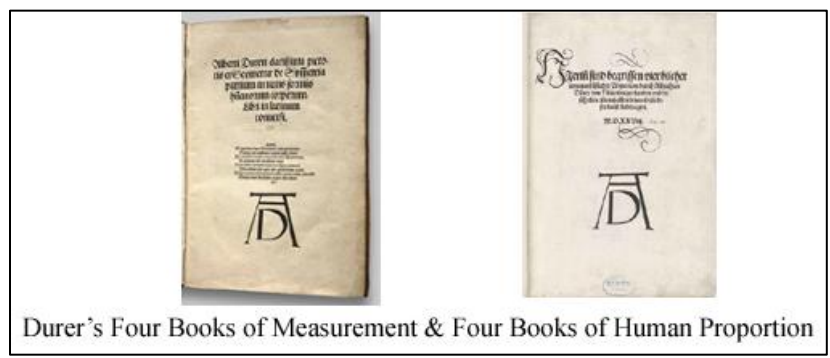

Durer's Four Books of Measurement & Four Books of Human Proportion

These books were written in German instead of the traditional Latin, so they were accessible to craftsmen, including printmakers, painters, goldsmiths, sculptors, and stonemasons. To do this Dürer often had to invent a vocabulary for some technical phrases which only existed in Latin. (As late as the 19$^{th}$ century, German mathematicians such as Carl Friedrich Gauss still wrote in Latin.)

The Books of Measurement were a textbook in practical geometry so that the artist could "*combine practical skills with theoretical skills*". These were among the first books on geometry and mathematics ever written for the general public. For his sources, Dürer had a copy of Euclid, but most of the material apparently came from workshops and pamphlets such as *Geometrica deutsch* which was printed in about 1486. It was preceded by *Geometria Culmensis* from about 1400. He also had access to *Fialenbuchlein* which was published in Nuremberg in 1486 and any of the books in the library of his wealthy friend Pirckheimer – who had a collection of books and manuscripts on science and mathematics.

**Four Books on Measurement**

In the introduction shown below, Dürer states that it is "*a manual of measurement of lines, areas, and solids by means of compass and ruler assembled by Albrecht Dürer for the use of lovers of art with appropriate illustrations to be printed in the year MDXXV*"

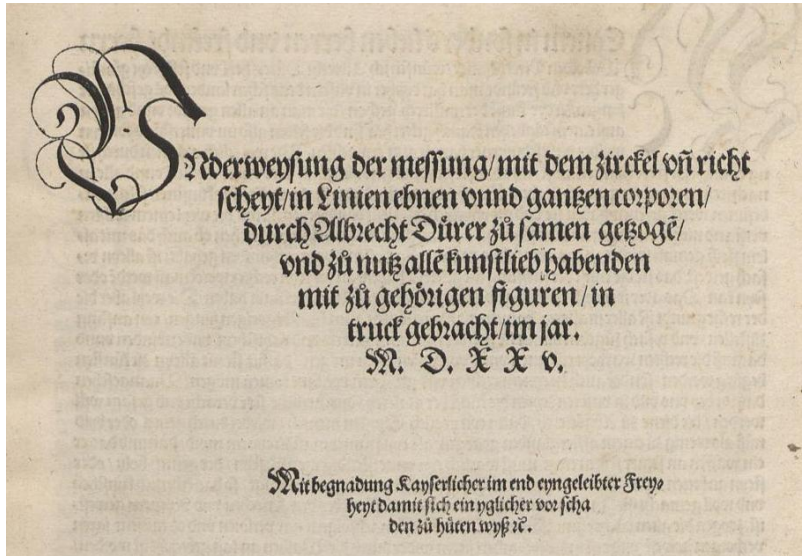

The books progress from "*Lines to Planes to Whole Bodies*". The first section covers conics, the spiral of Archimedes and the conchoid. The second and third sections deal with planar constructions and their applications, and the fourth section deals with Platonic and Archimedean solids. Below is a (rotated) fig. 17 from Book 1 showing a sinusoidal curve which Dürer derives from 12 points on a circle. Dürer calls this spiral a 'screw-line'. It was apparently used by stonemasons. The text below says "*This is the spiral and its ground plan as described above*".

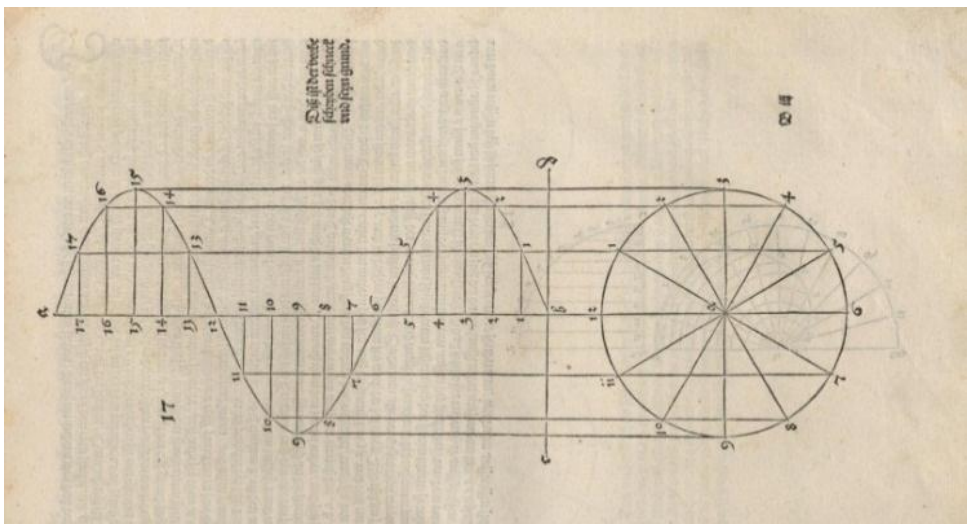

**Part II – The polygon constructions from the second book of measurement**

The second book contains illustrations and directions for the construction of geometrical objects such as 'regular' polygons with compass and (unmarked) straightedge. Dürer was very modest and never claimed that his constructions were new or unique – although some appear to be unique. Most of the exact constructions were passed down from Euclid and the approximate constructions were typically those used in the workshops of engravers, cabinetmakers or masons.

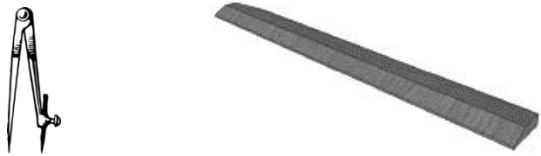

Compass and straightedge constructions data back to Euclid of Alexandria who was born in about 300 B.C. In terms of regular polygons, the Greeks could construct polygons with $2^k$ sides for k >1. They could also construct the regular triangle and the regular pentagon but these were the only 'prime 'polygons that they could construct. They knew that the construction of a regular n-gon could be combined with the construction of a regular m-gon to yield a regular nm-gon when n and m were relatively prime. Dürer shows how to construct a 15-gon using a triangle and pentagon.

Regular polygons had been used throughout the Middle Ages in architecture and decoration as well as in the design of cities and fortifications. Dürer knew that there were a limited number of 'exact' constructions and in the case of the pentagon he gives both an exact and an approximate construction. The approximate construction, based on a single compass setting, may be due to Pappus of Alexandria (c. 290 – 350 A.D.). Dürer preferred such constructions when possible. It yields an equilateral pentagon which Dürer apparently thought was also equiangular.

In the Second Book, Dürer covers all the 'regular' polygons from the triangle to the 16-gon. (He leaves it to the reader to get the dodecagon by bisecting the arcs of the hexagon.) Dürer also mentions how to construct the 28-gon. He was probably aware that the St. Lorenz church in Nuremberg had a beautiful rosette window with stonework based on a regular 28-gon. This window was apparently constructed in about 1470. The church and window are shown below. Note the transition from an inner octagon to the outer 28-gon.

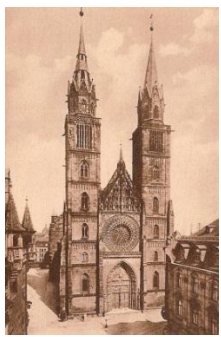 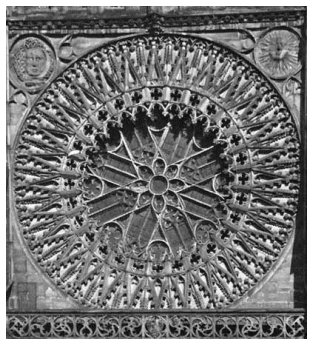

Below is a list of the polygon constructions in the order they are given in the Second Book.

| Number of sides | Dürer's figure | Regular (Y/N) | Thumbnail of figures |
|---|---|---|---|
| 6 | 9 | Y | 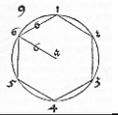 |
| 3 | 10 | Y | 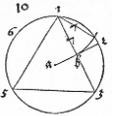 |
| 7 | 11 | N | 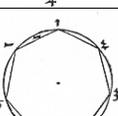 |
| 14 & 28 | 12 | N | 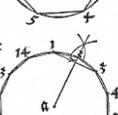 |
| 4 | 13 | Y | 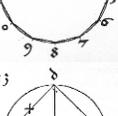 |
| 8 & 16 | 14 | Y | 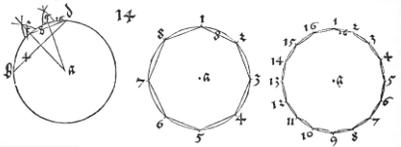 |
| 5 & 10 | 15 | Y | 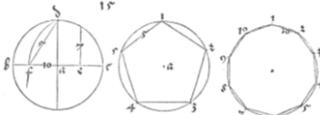 |
| 5 | 16 | N | 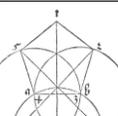 |
| 15 | 17 | N | 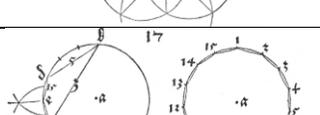 |
| 9 | 18 | N | 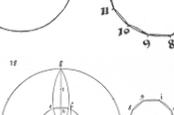 |
| 11 & 13 | 19 | N | 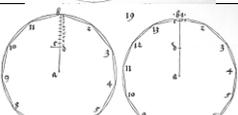 |

For each construction we will give Dürer's original 1525 text and figures along with a translation. A photocopy of the original book is available on-line at the Saxon State library in Dresden. The most notable mathematical studies were by H. Steigmuller in 1891 and Max Steck in 1948. The definitive 1943 study by Edwin Panofsky focuses mainly on his art. In 1977, W.L. Strauss published a copy of the Four Books with a side-by-side English translation and some analysis based on Steigmuller. More recent analysis and partial translations have been published by Donald Crowe and Dan Pedoe.

A copy of the 1525 text is also available at wikisource and Google will attempt a translation, but the resulting translation of Dürer's medieval German to English is almost comical. The passages needed here are relatively short and the accompanying diagrams make translations fairly routine, but there are places where Dürer's text and his diagrams seem to be at odds. In these cases we will present views and opinions from as many sources as possible.

After presenting the text and translations, we will carry out the non-regular constructions with Mathematica and in a related paper, we will perform an analysis of these non-regular cases using the Tangent map (outer billiards map). We encourage the reader to download the Mathematica code from DynamicsOfPolygons.org and perform their own analysis.

Below is a page from the original manuscript showing figures 9,10,11 and 12. Dürer presented the polygons in logical order – starting with the hexagon and using the odd vertices to get a regular triangle and then bisecting a side of this triangle to construct the heptagon and finally bisecting the arcs of the heptagon to get the 14-gon.

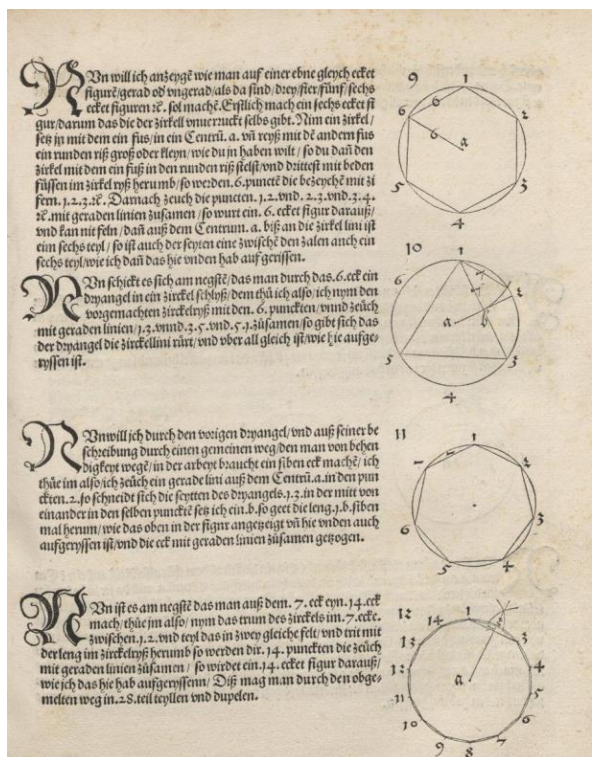

For the most part Dürer knew which constructions were approximate and which were exact, but he was more concerned with the practical aspects of the constructions. Whenever possible he gave constructions which used a fixed compass setting. Interest in such 'rusty compass' constructions was fueled by the popularity of his work, but it was later shown that any construction with variable compass can also be carried out with fixed compass. Some of the early compasses did not even have a hinge but it is still possible to transcribe a given circle at another location, just as if there was a hinge.

Approximate constructions were considered to be 'mechanical' and not demonstrative. Dürer noted that 7, 11 and 13 were in this category. He does not state explicitly that the 9-gon construction is mechanical but it is based on 'workshop' methods which were well known to be approximate. We remarked earlier that Dürer did not state that the fixed compass pentagon construction was approximate.

◆**The Regular Hexagon – Fig. 9**

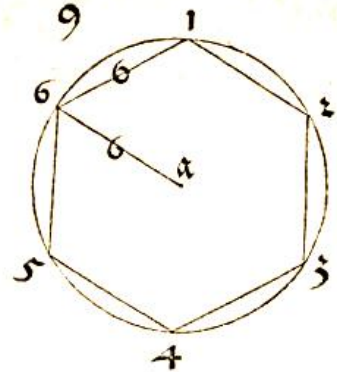

*"Now I will explain how to construct figures with an even and uneven number of corners, that is three, four, five and six-cornered figures, etc. First I shall construct a figure with six corners because it can be made without changing the opening of the compass. Take a compass and place one leg on center point a and with the other side draw a circle of the desired size, large or small. Then move the compass without changing its opening, to the periphery of the circle, and mark off this distance on it. Then continue along the periphery in the same manner. You will then obtain six points which you will mark with numerals 1,2,3,etc. Connect points 1 and 2, points 2 and 3, points 3 and 4,etc. with straight lines. You will obtain a six-cornered figure. Therefore the line from the center a to the periphery is one-sixth, and also each line between the numerals is one-sixth."*

◆ **The Regular Triangle – fig. 10**

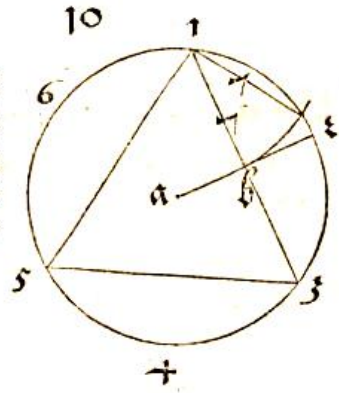

*"Now I will use this six-cornered figure to construct a triangle within a circle. I use the six points of the preceding figure, but connect points 1 and 3, points 3 and 5, and points 5 and 1 with straight lines. This will result in a triangle that touches the periphery and is symmetrical in every respect and position, as shown in the diagram."*

**Commen**t: Dürer combined the triangle and heptagon constructions in fig. 10 and showed the results of the heptagon construction in fig. 11 below.

◆ **The Regular Heptagon (approximate)- fig. 11**

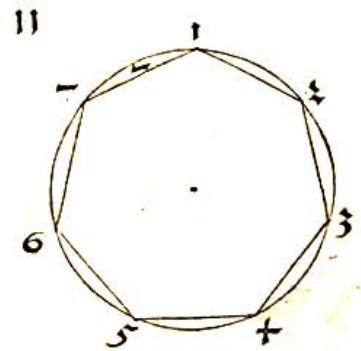

*"Now I shall show a simple way of using the triangle of the preceding figure to construct a seven-cornered figure. I draw a straight line from the center point a to point 2 which cuts the side 1,3 of the triangle in half. Where this occurs I mark b. The length 1b will then fit seven times around the periphery, as shown in the figure."*

**Comment**: It had usually been assumed that Dürer obtained this construction from the pamphlet *Geometria deutsch* but Donald Crowe points out that the construction there is different and does not conform to the strict compass and straightedge criteria, so this construction may be Dürer's own. Indeed Kepler refers to this as 'Dürer's Rule'. But according to J. Tropfke this construction has a long history and was known to Heron of Alexandria.

◆ **The Regular 14-gon (approximate)- fig. 12**

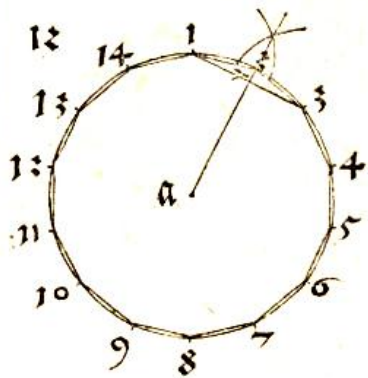

 *"Next I will show how to construct a figure with fourteen sides from only seven sides. By means of a circle, divide one of the sides in half. Then mark off this length along the periphery. You will have fourteen points which when connected will form a figure with fourteen corners as shown in the diagram. One can use the same method to construct a figure with twenty-eight sides or twice as many."*

**Comment**: There has been some contention among translators about Dürer's meaning when he says to "mark off this length along the periphery". There is little doubt that Dürer intended for the reader to use the bisected arc rather than the bisected side. He certainly knew the difference and the diagram seems to make this clear. The procedure for bisecting a cord and an arc are the same. It's just a matter of what point you chose after drawing the two 'vesica pices' arcs as shown above. It is natural to choose the point that is already on the original arc, and that is clearly what Durer did here. The only exception is in figures 10 and 11 where Durer deliberately chose the length rather than the arc, so that the resulting figure has 7 sides instead of 6.

In the next two constructions Dürer leaves no doubt that to double the number of sides, arcs must be bisected. Both H. Steigmuller and W.L. Strauss seemed to be of the opinion that Dürer meant for the reader to use the given length, but Strauss was not a mathematician and his analysis was typically based on that of Steigmuller. Donald Crowe remarks that Steigmuller apparently did not have his own copy of the manuscript so he had to rely more on the text (which he could copy from the library ) rather than the diagrams. Dürer of course assumed that the reader has the diagram as well as the text and is using common sense. As we will see with the 13-gon, Dürer sometimes assumes that the diagrams can speak for themselves.

◆**The Regular Quadrilateral – fig. 13**

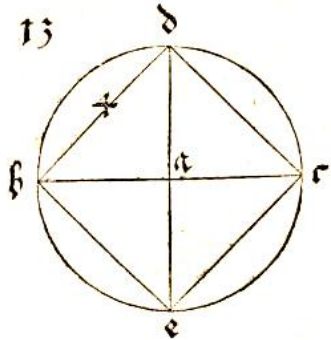

"Now I will construct a quadrangle with a compass. First draw a circle. Through its center a draw a horizontal line. Where this line crosses the periphery I mark the points b and c. Then I draw a vertical line from center a up and down to the periphery, at right angles to the horizontal line, and I mark the point where it crosses the periphery. On top I mark d and on the bottom e. Then I connect points b and d, d and c, c and e and c and b by straight lines. They will form a square symmetric with the circle, as shown in the diagram below."

◆**The Regular Octagon – fig. 14**

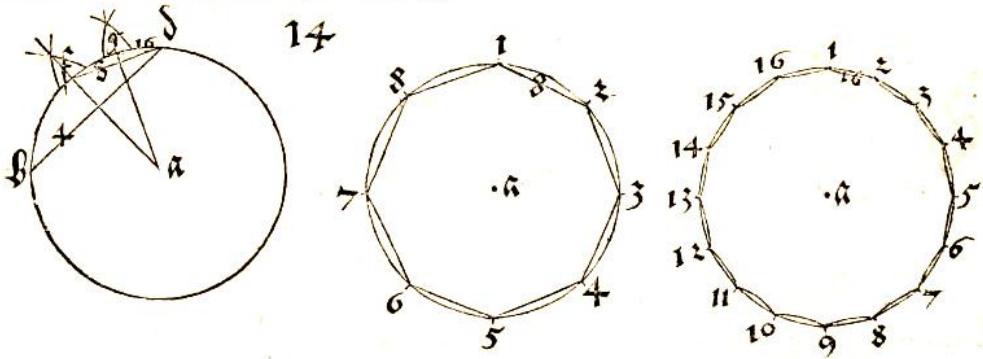

*"Now it is proper to show how to construct an eight-cornered figure. Using the preceding figure, retain the side bd and divide the arc above it into two halves. Mark the midpoint f. Then connect f and d and it will give you one side of the eight cornered figure. To obtain a sixteen cornered figure, divide the arc fd into two halves and mark the midpoint g. The line connecting g and d will represent one side of the sixteen sided figure. This is demonstrated in the three figures below."*

◆**The Regular Pentagon- fig. 15**

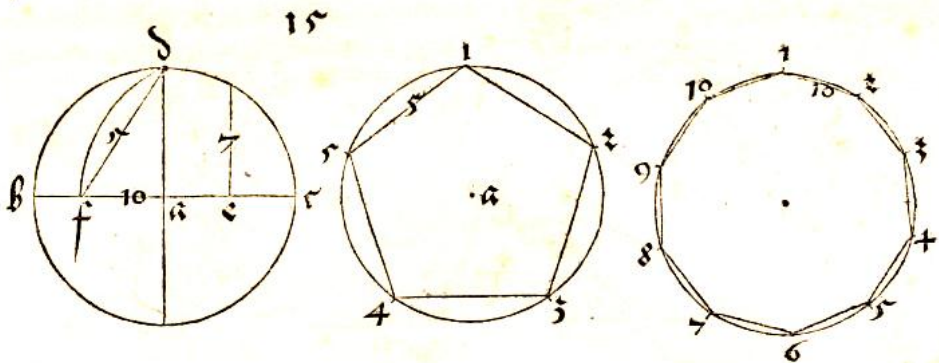

*"Now it is necessary to construct a pentagon within a circle. To do this, first draw a circle with center a and draw a horizontal line through this center. Where this horizontal line crosses the periphery of the circle mark points b and c. Then draw a vertical line through the center a at right angles to the horizontal line. Where it crosses the periphery at top mar the point d. Now draw a straight line ed and place one leg of a compass on e and the other on d and draw an arc down to the horizontal line bc. Where the arc crosses the line mark point f. Then connect f and d. The line fd represents one edge of the pentagon, whereas line fa is equal to one side of a ten-sided figure. Then divide line ac into two halves and erect a vertical line up to the periphery. Its length is equal to approximately one-seventh of the circle. This is shown in the following diagram."*

**Comment**: This is an old construction credited to Caludius Ptolemy (c.80-198). In the first of these three figures the '5', '10' and '7' show that this one diagram contains an (exact) side of a pentagon and decagon and an (approximate) side of a heptagon. The heptagon side given here is $\sqrt{3}/2$ just as in Dürer's method, and he notes that it is approximate.

Dürer of course intended that the point *e* is the bisector line *ac*. There is much speculation as to why Dürer chose this construction rather than Euclid's construction which uses the 'golden-ratio' proportion. The speculation stems from the fact that Dürer makes no mention of the golden ratio, although he was no doubt aware of its use in Italian art. It may be that Dürer simply did not feel comfortable with the precepts of the 'divine' ratio. The German architects had their own 'divine' ratio which was the *vesica pices* ratio of $1:\sqrt{3}$ (more on this below). In addition Ptolemy's construction is simpler than Euclid's and these constructions were just a preliminary step in his program.

◆**The Regular Pentagon (approximate)- fig. 16**

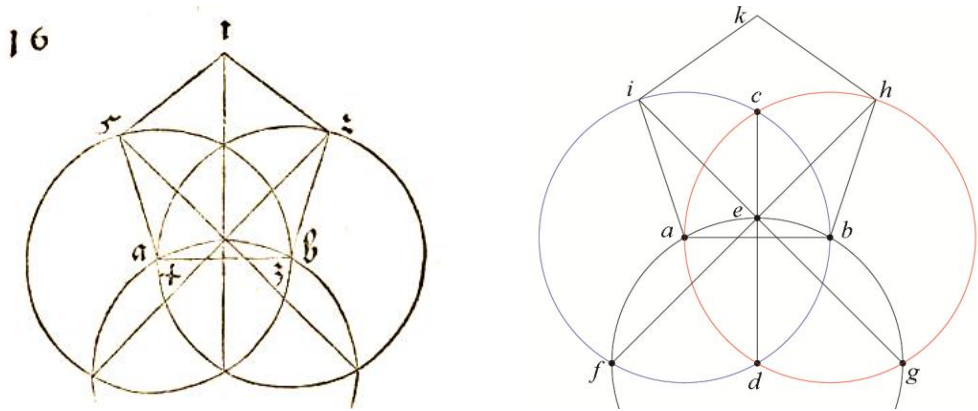

   "*Construction of a pentagon without changing the opening of the compass is accomplished as follows. Draw two circles which overlap so that the periphery of each touches the other's center. Then connect the two center points a and b with a straight line. The length of this line is equal to one side of the pentagon. Where the two circles cross mark c at top and d at the bottom and draw the line cd. Then take the compass without changing the opening and place one leg on d and with the other draw an arc through the two circles and their centers a and b. And where the periphery is crossed by this arc mark points e and f. Where the vertical line cd is crossed mark the point g. Then draw a line eg and extend it to the periphery of the circle. Mark that point i. Then connect i with a and h with b and it will give you three sides of the pentagon. Then erect two inclined lines of equal length from i and h until they meet at top. You will have pentagon as I have drawn here.*"

**Comment:** Points *c* through *i* are not labeled in Dürer's diagram so we have included a diagram on the right above showing the locations of these points. This construction and diagram are taken directly from *Geometria Deutsch,* but in that pamphlet the reader is told in the last step to draw a circle at *h* and take *k* to be the point where this circle intersects the line of symmetry *cd*. Dürer apparently felt the need to instead draw two symmetric circles at *h* and *i*. These last two circles are not drawn above.

This fixed compass construction yields an equilateral pentagon which is not equiangular. Dürer did not state that it was approximate and he may have believed that the angles were equal. This construction apparently dates back to Pappas of Alexandria. Because Dürer's works were widely distributed, many of the constructions are credited to him and Donald Crowe points out that P. Cataldi in 1620 derived extensive equations to show that 'Dürer's pentagon' is not regular.

◆**The Regular 15-gon (approximate or exact depending on which pentagon is used)- fig. 17**

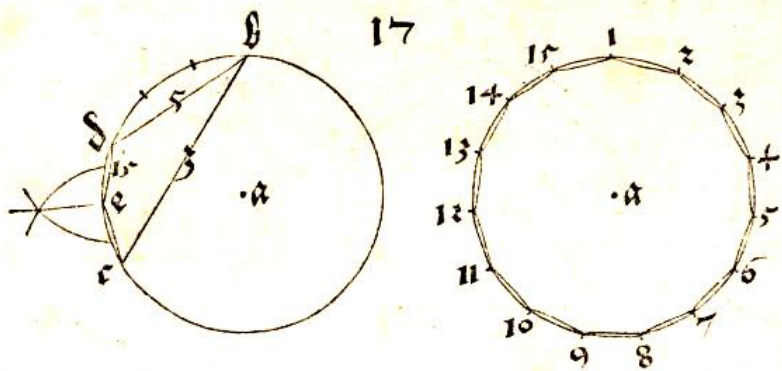

  "By using the construction of a triangle, as outlined above, it is possible to convert this pentagon into a 15-sided figure. First draw a circle with center a. then draw one side of a the triangle, marked b at top and c at bottom. Then take the length of one side of the pentagon and place one end of it on point b and the other on the periphery of the circle. Where it touches the periphery mark the point d. Divide the remaining distance between d and c be means of a compass into two halves and mark the midpoint e. If you then connect points e and c, it will give you the length of one side of the fifteen-sided figure, as shown in the following diagram."

**Comment**: This construction was described in Euclid. Note that if the second edge of the pentagon is drawn, then the arc bisection can be avoided, but neither Euclid nor Dürer do this.

◆**The Regular 9-gon (approximate)- fig. 18**

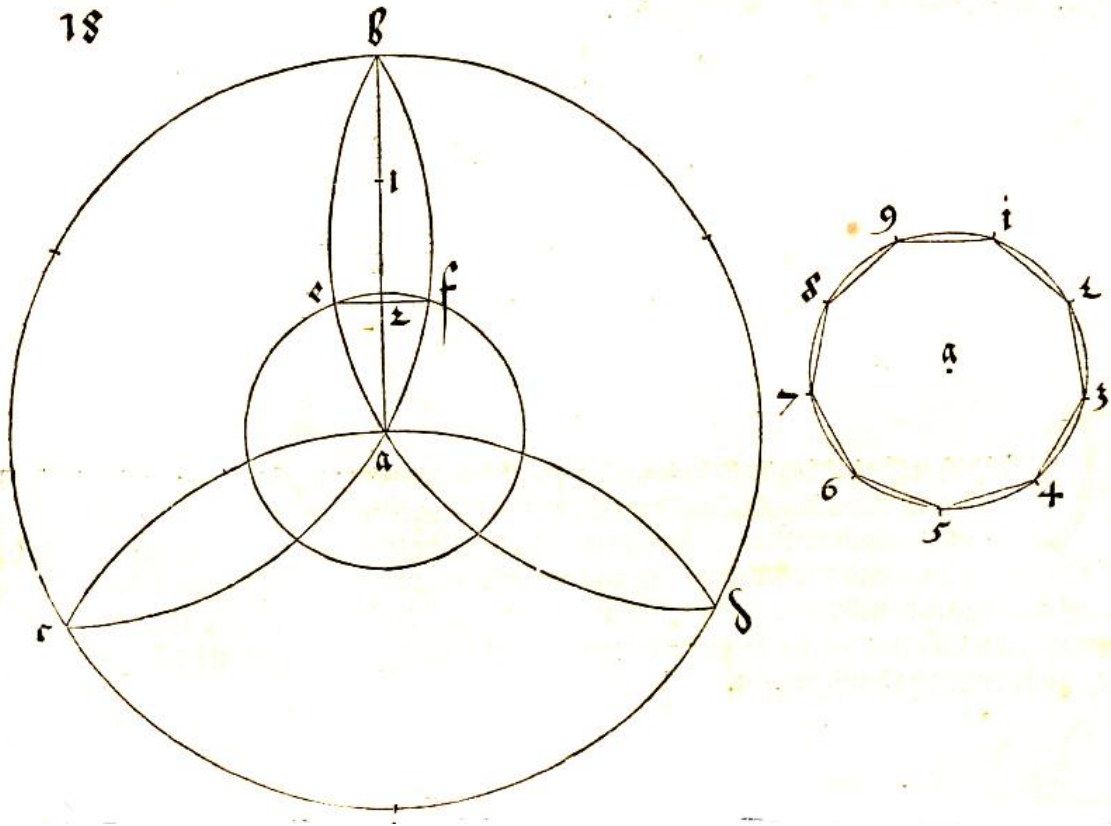

In neün eck ist durch ein dryangel zů finden/also/ Reiß auß einem Centrum. a. ein grosse zirckellini/ darein reiß mit vnuerrucktem zirckel/drey fischs blosen/der obern ende an der zirckellini sey. b. der andern end auf den seyten sey.c. d. Darnach reiß in der obern fischblosen/ ein aufrechte gerade lini.b.a. dise lini teil mit zweyen punckten.1.2. in drey gleiche felt/also das 2. der negst punckt beym.a. sey/ vnnd far durch den punckten.2. mit einer geraden zwerch lini zů gleichen wincklen.b.a. vnd wo sie die blosen lini zůbeden seyten durchschneidet/da setz.e.f. Darnach nym ein zirckel/setz jn mit dem ein fůß/in das Centrū.a. vnd den andern in den punckten. e. vnd reiß durch das.f. zů ring herumb / ein zirckellini so geet die leng.e.f. zů neün mal in disem zirckelriß herum/solchs hab jch hernach aufgeryssen.

*"You can construct a nine-sided figure based on a triangle. Draw a large circle with center a. Then without changing the opening of the compass draw three 'fish bladders' whose upper end on the periphery you will mark b. Mark the other c and d. Within the upper fish bladder draw a vertical line ba and divide this line with two points 1 and 2 into three equal parts. Point 2 should be closest to a. Then draw a horizontal line through point 2 at right angles to the vertical line ba. Where the horizontal line crosses the fish bladder mark points e and f. Then place one leg of the compass on center a and the other on point e and draw a circle through f. Line ef will then represent one of nine sides which will compose a nonagon inside this smaller circle, as shown in the diagram below."*

**Comment**: Dürer's *fischblosen* (fish bladders) are known in Latin as *vesica pices* and in Italian as *mandorla* (almonds). They occur quite naturally when bisecting a given length or arc with a compass. Vesica pices have appeared in art and architecture dating back to early Egypt. Dürer was familiar with the work of the Roman architect Marcus Vitruvius (c. 80 BC- 15 BC) who formulated a science of proportion based on geometric form and symmetry in the spirit of Pythagoras. Vitruvius compiled a set of Ten Books on Architecture and he also studied human proportion. The ratio of the length to width of the vesica pices is $\sqrt{3}$ and this proportion was common in medieval and Gothic architecture. In Cesar Casariano's *Vitruvius* of 1521 he calls this the 'rule of the German architects'.

In Dürer's unpublished notebooks at the British Museum there are drawings which suggest that he did not simply reproduce the given *ef* interval around the inner circle, but instead he used all three fish bladders and obtained the remaining vertices by bisection. See the details below in the analysis section.

Dürer did not state explicitly that this construction is approximate, but most likely he knew it was not 'demonstrative'. Edwin Panofsy notes that " ..*his approximate construction of the enneagon is not described in any written source but was taken over directly, as we happen to know, from the tradition of the 'ordinary workmen'* "

◆**The Regular 11-gon and 13-gon (both approximate) – fig. 19**

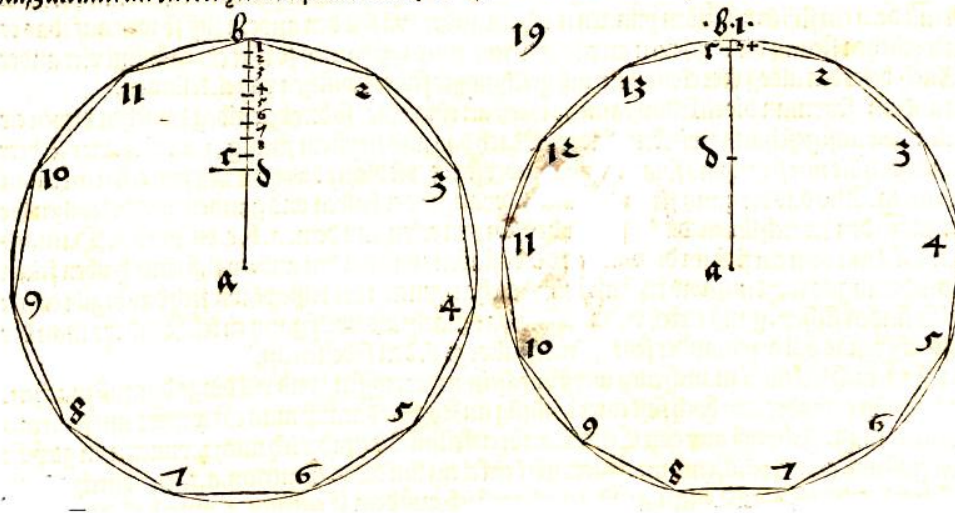

*"To construct an eleven-sided figure by means of a compass, I take a quarter of a circle's diameter, extend it by one-eighth of its length and use this for construction of the eleven-sided figure. This is a mechanical not a demonstrative construction."*

*"If a thirteen-sided figure is required, I first draw a circle with center a. Then I draw the radius ab and cut it in half at d. I use the length cd to mark off thirteen pieces along the periphery of the circle. But this method is mechanical and not demonstrative."*

**Comment**: Dürer only devoted a few sentences to both of these but the resulting constructions are surprisingly accurate. There is no known historical precedent for either, and they may be unique to Dürer. The 11-gon construction is more accurate than the traditional method which is based on 'Renaldini's Rule'.

Both Steigmuller and Strauss seemed to be confused by Dürer's intentions for the 13-gon and they omitted any reference to the 'length cd'. Dürer clearly states to use **'leng .cd'** but both Steigmuller and Strauss say "use this length" which of course implies that the reader should simply cut *ab* in half and use that length. Stegimuller knew that this would create a large error and he remarked "*Even in Dürer's figure the mistake is disturbingly apparent at first glance."* Strauss repeats this opinion when he discusses Steigmuller's analysis: "*… for the regular thirteen sided polygon the equivalents are: actual .478635; Dürer .5 but this difference makes the method imperfect as can be seen even in Dürer's figure."*

In fact the picture is very accurate and Strauss was just repeating the opinion of Steigmuller. Indeed if the side was .5, the figure would look very weird because the 13$^{th}$ side would be about .2184. This is clearly not what Dürer intended or what is drawn above. There is no doubt that *bd* is intended to be ½ and *cd* is smaller. Why else would there even be a point marked *c* ? We will address this issue in the analysis section below.

**Part III - Analysis of the approximate constructions**

◆**Dürer's Heptagon**

Dürer's construction of the heptagon has 6 equal sides but the base is slightly longer so the resulting polygon is not quite equilateral. This means that the corresponding 14-gon and 28-gon will suffer from the same imperfections since these are obtained by bisecting the arcs of the heptagon. It is impossible to construct a regular heptagon with compass and ruler. For a demonstration of Dürer's construction see the CDF manipulate at DürersHeptagon. The first stage of the manipulate is shown below next to Dürer's drawing which instructs the user to bisect the side of the equilateral triangle and use the compass to reproduce this length as an edge of the heptagon. Bisecting the arc rather than the side would lead back to the hexagon, so the difference between bisecting the side rather than the arc gives the heptagon.

Assuming that the radius of the circle is 1, the sides of the equilateral triangle are each $\sqrt{3}$ so the six sides will have length $\sqrt{3}/2 \approx 0.866025403784439$ which is within 3 decimal places of the actual which is $\approx 0.8677674782351162409$. The base will be slightly longer at $\dfrac{9\sqrt{39}}{64} \approx$ 0.878202843524774. The interior angles (other than the base) are $2\text{ArcSine}[\sqrt{3}/4] \approx 51° \, 19'$ 4.12517" compared with $360/7 \approx 51° \, 25' \, 42.85714"$

The exact coordinates of the vertices are: Mc=

$$\{\{0,1\},\{\tfrac{\sqrt{39}}{8},\tfrac{5}{8}\},\{\tfrac{5\sqrt{39}}{32},-\tfrac{7}{32}\},\{\tfrac{9\sqrt{39}}{128},-\tfrac{115}{128}\},\{-\tfrac{9\sqrt{39}}{128},-\tfrac{115}{128}\},\{-\tfrac{5\sqrt{39}}{32},-\tfrac{7}{32}\},\{-\tfrac{\sqrt{39}}{8},\tfrac{5}{8}\}\}$$

The vertices to 25-decimal places:

N[Mc,25]= {{0,1},{0.7806247497997997757308616,0.6250000000000000000000000},
{0.9757809372497497196635771,-.2187500000000000000000000},
{0.4391014217623873738486097,-0.8984375000000000000000000},
{-0.4391014217623873738486097,-0.8984375000000000000000000},
{-0.9757809372497497196635771,-.2187500000000000000000000},
{-0.7806247497997997757308616,0.6250000000000000000000000}}

### ⬥Dürer's Equilateral Pentagon

Below are two compass and straightedge constructions of equilateral pentagons. The one on the left is far from regular but the one on the right is 'almost' regular. This is the construction used by Dürer. It was well-known in Dürer's time and was apparently part of Greek folklore. Both constructions are based on the same principle – a chain of circles with constant radius. The radius determines the side of the pentagon. Dürer prefered constructions which did not involve changing the compass settings and these two constructions are both in this category of 'rusty compass' constructions.

In the example on the left, start with an arbitrary circle c1 and pick any point on c1 for the center of c2. Then choose an obtuse angle α ,and continue in this fashion until one of the circles intersects c1. This gives a path home. A similar construction works for any number of sides. The pentagon shown here is completely determined by the angles α and β.

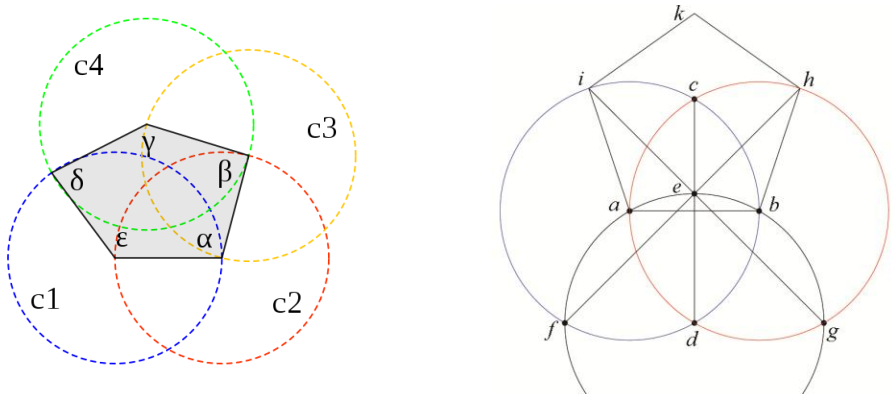

The only difference in Dürer's construction is that he provides a method for choosing a 'good' angle α. (His construction yields α ≈ 108°21'58.03259 compared to the exact 108°.) As on the left, start with 'vesica pices' circles centered at *a* and *b*. If *a* = {0,0}, and the radius is 1, then *c* and *d* are $\{1/2, \pm\sqrt{3}/2\}$. The lower circle at *d* determines the points *e*, *f* and *g* which are vertical and horizontal displacements of *d*. The lines *fe* and *ge* determine *h* and *i*. The final point *k* must lie on the extension of *cd*, so a circle centered at *h* or *i* will locate *k*.

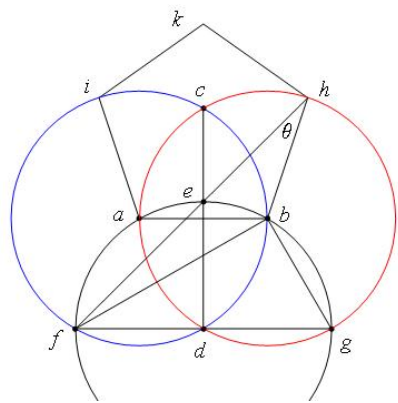

To find α = *abh*, note that the triangle *bfg* is a right triangle with base $\sqrt{3}$. This implies that in the adjacent triangle *bfh*, the lower angle is 15° and the Law of Sines gives the upper angle as $\theta$ = ArcSin[$\sqrt{3}$ Sin[15°]]. Therefore α = 180° -15° - θ -30° ≈ 108°21'58.03259 as indicated above. To find β = *bhk*, use the fact that triangle *hbc* is isosceles: β ≈ 107° 2'16.17332. This yields γ ≈ 109° 11'31.58819. The coordinates of *k* can be found using the isosceles triangle *ikh*.

The vertices of Durer's pentagon are given below in terms of $\theta = \text{ArcSin}[\frac{1}{2}\sqrt{\frac{3}{2}}(-1+\sqrt{3})]$

$a = \{0,0\}, b = \{1,0\}, h = \{ 1+\text{Cos}[\frac{\pi}{4}+\theta], \text{Sin}[\frac{\pi}{4}+\theta]\}, i = \{-\text{Cos}[\frac{\pi}{4}+\theta], \text{Sin}[\frac{\pi}{4}+\theta]\},$

$k = \{1/2, \text{Sin}[\frac{\pi}{4}+\theta]+\sqrt{1-(\text{Cos}[\frac{\pi}{4}+\theta]+1/2)^2}\}$

We list the vertices below to 25 decimal places, clockwise starting with k:

Mc = {{0.5000000000000000000000000, 1.528399804246314681857 0889},
{1.3150878974491474905777 1507449, 0.9490624936647088438139 91903739},
{1, 0}, {0, 0}, {-1.3150878974491474905777 1507449, 0.9490624936647088438139 91903739}}

Our convention is to have vertex 1 at {0,1} and a radius of 1, so we will scale Mc and then translate.

**s0 = SideFromRadius[1, 5]**= 1.17557050458494625833741190928 *side length for regular pentagon with radius 1*)  **Ms = s0*Mc; v = {0, 1} - Ms[[1]]; tr = TranslationTransform[v];**

**Mc = tr[Ms]** = {{0, 1},{0.958193290885377245697767, 0.318948145474695899702365},
{0.587785252292473129168706, -0.796741729085373237430852},
{-0.587785252292473129168706, -0.796741729085373237430852},
{-0.958193290885377245697767, 0.318948145474695899702365}}

### ◆Dürer's 9-gon

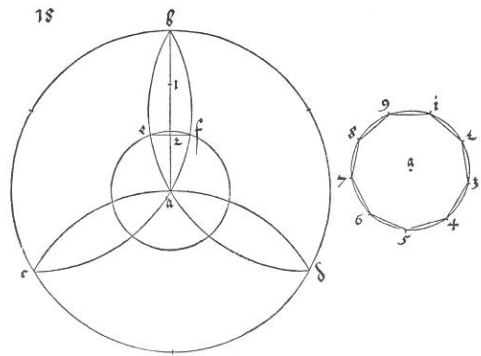

Dürer's drawing has three small marks where the compass should be placed to get the fish bladders- these marks can be seen above on the outer circle at 30°, 150° and 270°. The three fish bladders determine 3 of the 9 edges, and as indicated earlier, there is evidence from Durer's sketchbooks that he found the remaining vertices by bisecting the remaining arcs. This would ensure symmetry and minimize the error propagation. We will do this below.

To find the points *f* and *e,* use the formula for the intersection of two circles of radii R and r assuming that the R circle is at the origin and the r circle is at {d,0} . See wolfram math world.

The two intersections are at {x,y} where $x = (d^2-r^2+R^2)/2d$ and $y^2 = [4d^2R^2-(d^2-r^2+R^2)]^2/4d^2$

In our example r = 1, R = 3 and d = 3 so the top point is {x,y} = {17/36, $\sqrt{35}$/6}.

We need to rotate this by π/2 to get **p1= RotationTransform[Pi/2][{x,y}]** = {-$\sqrt{35}$/6, -1/6}. This is one of the vertices and two symmetric vertices can be obtained by rotations of ±120° , so

**f = RotationTransform[-2*Pi/3][p1]** = $\{\frac{1}{12}(-\sqrt{3}+\sqrt{35}), \frac{1}{12}(1+\sqrt{105})\}$ and *e* is a reflection, so

$e = \{\frac{1}{12}(\sqrt{3}-\sqrt{35}), \frac{1}{12}(1+\sqrt{105})\}$. Rotating *e* and *f* by 120° yields a total of 6 vertices as shown below. We will follow Durer and bisect the remaining three arcs. This can be accomplished by rotating {0,-1} by ±120°

This gives the following vertex list, where we begin with {0,-1} by convention and proceed clockwise:

$$Mc = \{\{0,-1\}, \{\frac{1}{12}(-\sqrt{3}-\sqrt{35}), \frac{1}{12}(1-\sqrt{105})\}, \{-\frac{\sqrt{35}}{6},-\frac{1}{6}\}, \{-\frac{\sqrt{3}}{2},\frac{1}{2}\}$$

$$\{\frac{1}{12}(\sqrt{3}-\sqrt{35}), \frac{1}{12}(1+\sqrt{105})\}, \{\frac{1}{12}(-\sqrt{3}+\sqrt{35}), \frac{1}{12}(1+\sqrt{105})\}, \{\frac{\sqrt{3}}{2},\frac{1}{2}\}$$

$$\{\frac{\sqrt{35}}{6},-\frac{1}{6}\}, \{\frac{1}{12}(\sqrt{3}+\sqrt{35}), \frac{1}{12}(1-\sqrt{105})\} \}$$

To put Mc in standard form, rotate about the x axis (and keep the clockwise orientation)

Mc = {{0,1},{0.637344215889041111113412312,0.770579230496331986017532},
{0.986013297183269340427880,0.166666666666666666666667},
{0.866025403784438646763723,-.500000000000000000000000},
{0.348669081294228229086568,-0.937245897163299865268419},
{-0.348669081294228229086568,-0.937245897163299865268419},
{-0.866025403784438646763723,-0.500000000000000000000000},
{-0.986013297183269340427880,0.166666666666666666666667},
{-0.637344215889041111113412312,0.770579230496331986017532}}

At this scale it is indistinguishable from the regular 9-gon but there are two sets of interior angles. The three filled regions above are the original sectors determined by e and f . Their angles are 40° 48' 42.70877" and the remaining 6 sectors are a little smaller at 39° 35' 38.64562".

There is a CDF manipulate that performs this construction by projecting the points to the outer circle. The manipulate does not use bisections, so it yields a slightly different result. It can be accessed at Wolfram's Mathematica site: Dürer9-gon.

It may have been tempting for Dürer to reproduce the angle of 40° by a trisection of 60°, but there is no evidence that he attempted to do this. It is easy to prove that exact trisection of 60° is not possible with compass and straightedge, because it must be true that $\cos 3\theta = 4\cos^3\theta - 3\cos\theta$ and when $x = 2\cos 20°$, this becomes $x^3 - 3x - 1 = 0$ which is irreducible. Since $\cos 20°$ is not constructible, 20° is not constructible. See the Appendix for Dürer's approximate trisection method. There we show that his tisection method applied to 60° would yield 20° to .01 seconds of arc which is pretty amazing.

◆**Dürer's 11-gon**

*"To construct an eleven-sided figure by means of a compass, I take a quarter of a circle's diameter, extend it by one-eighth of its length and use this for construction of the eleven-sided figure."*

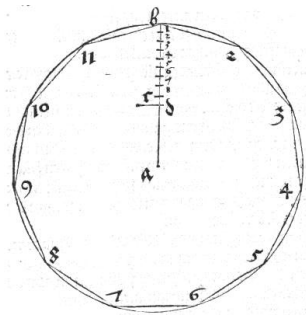

Using radius 1, the side length is $1/2 + (1/8)(1/2) = .5625$ which is a good approximation to the actual which is $\approx 0.5634651136828$. The corresponding interior angle is $\theta = 2\text{ArcSin}[(9/32)] \approx 32° 40'10.72402$ compared to $360/11 \approx 32°43'38.18182$. With 10 sides of length 9/16, the base will be slightly longer at $0.5731083727308245$

**theta = 2*ArcSin[(9/32)]; M1 = Table[RotationTransform[-k*theta][{0, 1}], {k, 0, 5}];**
**M2 = Reverse[Table[RotationTransform[k*theta][{0, 1}], {k, 1, 5}]]; Mc = Join[M1, M2];**

N[Mc,25] = {{0, 1}, {0.5397944249806905656719474,0.8417968750000000000000000},
{0.9087945201823345070492552,0.4172439575195312500000000},
{0.9902463492325366708055095,-0.1393275558948516845703125},
{0.7583780443458815285647082,-0.6518149598268792033195496},
{0.2865541863654123089264796,-0.9580640366261832241434604},
{-0.2865541863654123089264796,-0.9580640366261832241434604},
{-0.7583780443458815285647082,-0.6518149598268792033195496},
{-0.9902463492325366708055095,-0.1393275558948516845703125},
{-0.9087945201823345070492552,0.4172439575195312500000000},
{-0.5397944249806905656719474,0.8417968750000000000000000}}

◆ **Dürer's 13-gon**

*"If a thirteen-sided figure is required, I first draw a circle with center a. Then I draw the radius ab and cut it in half at d. I use the length cd to mark off thirteen pieces along the periphery of the circle. But this method is mechanical and not demonstrative."*

Dürer combined the 11-gon and 13-gon in Figure 19 as shown here.

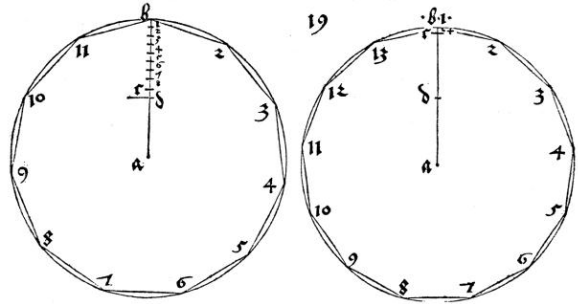

Concerning the 13-gon, we know that Dürer wants the reader to use the length cd which is slightly less than bd (which we will assume is ½) The question is, how much smaller? The answer to this question seems to depend on how we interpret the symbols enlarged below. It looks like a '2' and a small '+' sign but the most natural reading is the number 24. It is often the case in Dürer's text that '4' resembles '+'. This is the interpretation suggested by K. Hunrath in his 1905 article in *Bibliotecha Mathematica*. In the words of Donald Crowe "*Hunrath's reading is by far the most persuasive*". So Hunrath suggested that Dürer intended to decrease *bd* by one half of 1/24 to give an edge length of 1/2 - 1/48 = .4791$\overline{6}$

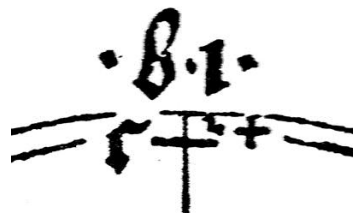

Another possibility has to do with the 11-gon on the left where Dürer had divided the interval *cb* into 8 equal segments and then added one of these 1/16 lengths to *cb* to get an edge of 9/16. Doing this one more time to get one half of 1/16, would give ½ - 1/32 = .46875. However the table below shows that this would be inferior to the 1/24 choice.

| (i) Edge length is 1/2 | 12 edges are .5 and edge 13 is about .2184 |
|---|---|
| (ii) Edge length is 1/2 - 1/32 = .46875 | 12 edges are .46875 and edge 13 is about .59617 |
| (iii) Edge length is 1/2 - 1/48 = .4791$\overline{6}$ | 12 edges are .47916 and edge 13 is about .4722 |
| Actual edge length ≈ .4786313285751 | |

So 1/24 is a very good choice and shows that Dürer took great pains to find accurate results. The real mystery is why Dürer was so cryptic about this construction. There are other figures which are missing complete descriptions so this omission is not totally unique. No changes to this text were made in the 1538 edition.

Since this construction is very accurate, it is unlikely that Dürer could test whether it was exact or not. This applies to the 11-gon as well. However Dürer probably had an intuition that there were no 'simple' constructions which would yield a perfect 11-gon or a perfect 13-gon and he labeled both of these accordingly.

Assuming that Durer wants the reader to reduce the half-radius *db* by 1/24 of *db*, the side length *cd* will be 23/48 with interior angle θ = 2ArcSin[23/96] ≈ 27° 43' 26.03764'' compared to the actual which is 360/13 ≈ 27° 41' 32.30769". With 12 sides of 23/48, the 13$^{th}$ side (the base) will be a little short at 0.47220444562087881310469856996.

**theta = 2*ArcSin[(23/96)]; M1 = Table[RotationTransform[-k*theta][{0, 1}], {k, 0, 6}];**
**M2 = Reverse[Table[RotationTransform[k*theta][{0, 1}], {k, 1, 6}]]; Mc = Join[M1, M2];**

N[Mc,25] = {{0,1},{0.4652113226503646256137798,0.8851996527777777777777778},
{0.8236098025567870259889791,0.5671568505557966820987654},
0.9929068998449193496505630,0.1188944415872806711288473},
{0.9342318834099777427717181,-0.3566662137352594175438313},
{0.6610565777719639892234987,-0.7503360586993132944627397},
{0.2361022228104394065523493,-0.9717282235192974088943263},
{-0.2361022228104394065523493,-0.9717282235192974088943263},
-0.6610565777719639892234987,-0.7503360586993132944627397},
{-0.9342318834099777427717181,-0.3566662137352594175438313},
{-0.9929068998449193496505630,0.1188944415872806711288473},
{-0.8236098025567870259889791,0.5671568505557966820987654},
{-0.4652113226503646256137798,0.8851996527777777777777778}}

**Appendix A: Dürer's trisection (approximate) – fig. 20**

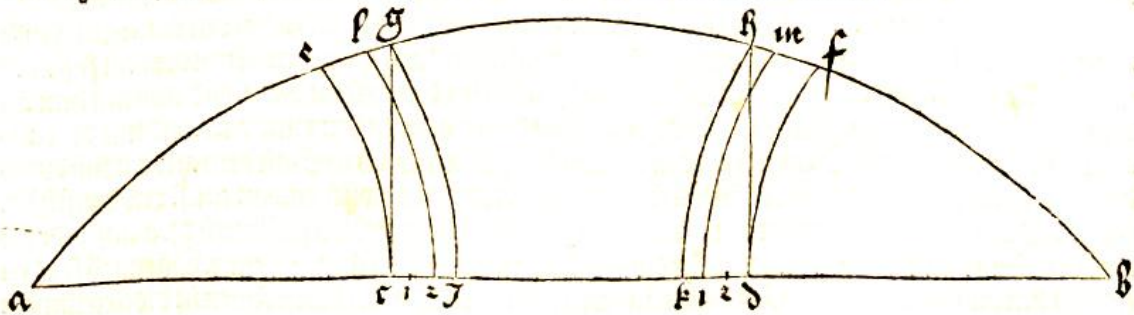

*" I divide an arc of a circle into three parts by the following method. The arc ab is connected by a straight line and, as I explained before, I divide this line into three equal parts with points c and d. Then I place one leg of the compass on a and with the other I draw an arc from c to the periphery. Where it crosses the periphery I mark the point e. Then I place one leg of the compass on b and draw an arc from d to the periphery. Where it crosses I mark the point f. Then I erect uprights from c and d to the periphery and where they cross I mark points g and h. The three arcs ae, gh and bf will be equal and this leaves just the two arcs eg and hf. I place one leg of the compass on a and the other on g and draw an arc down to ab, There I mark the point i. Then I place one leg of the compass on b and draw an arc down to ab and mark the point k. Then I divide ci and kd into three equal parts. Now I place one leg of the compass on a and the other on the point next to i and draw an arc up to the periphery and mark this point l. Then I place on leg of the compass on b and the other on the point next to k and again draw an arc up to the periphery and mark the point m. I have thus divided the large arc into three equal parts at points l and m as shown in the diagram. Whoever wants it more exact is welcome to look for a proof."*

Dürer knew that Euclid did not describe any method for trisecting an angle so he was aware that there was no known procedure that would give exact results. This probably led him to believe that his construction is not exact, but he did not venture an opinion in the text - except for his 'disclaimer' in the last sentence. For some reason this last sentence was omitted by Strauss who did remark that the method was approximate – yet highly accurate.

This construction appears to be unique to Dürer, who may have felt that there should be a connection between trisecting a cord and trisecting an angle. For small angles they are almost the same. The Greeks were adept at trisecting lengths and typically it takes 2 circles and 4 lines. What makes this angle trisection scheme so clever is that Dürer does it recursively with two successive trisections of intervals. This yield a 'second' order' correction to the first trisection. The end result is a very accurate trisection – even for larger angles.

**Example**: It is only necessary to do one of the two trisections. Below is a diagram showing the minimum steps needed to carry out Dürer's trisection of the angle AOB. (In this example AOB is 100°). The labeling here is the same as Dürer's, where we use the point X for the 'point to the right of k". This is the final point in the construction. Dürer extends this point up to the arc at m but that extension is not necessary if we only want to trisect the angle.

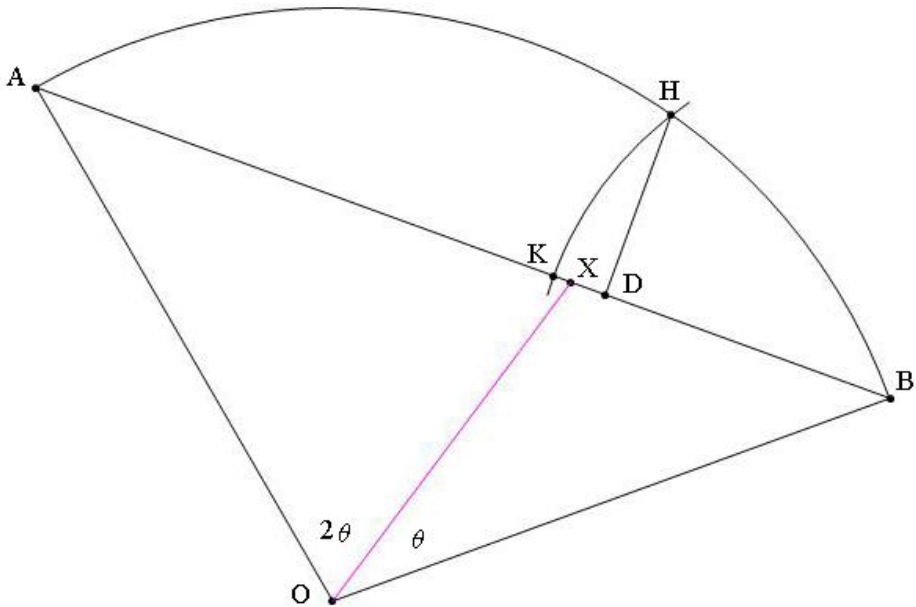

Step 1: Trisect the length AB to obtain the point D
Step 2: Construct a perpendicular at D which meets the arc at H
Step 3: Place the point of the compass on B and draw an arc from H to obtain the point K
Step 4: Trisect the length KD to get X

To do this with Mathematica, the only non-trivial step is finding H. Here are the calculations:

(*Points A and B*) **p1 = {Cos[20 Degree], Sin[20 Degree]}; p2 = {Cos[120 Degree], Sin[120 Degree]};**
**p3 = {0,0}; slope1 = (p1[[2]] - p2[[2]])/(p1[[1]] - p2[[1]]); d1 = EuclideanDistance[p1, p2]**
**Deltaxy[d_, slope_] := {d*Cos[ArcTan[slope]], d*Sin[ArcTan[slope]]};** (*to extend along a line*)
**p4 = p1 - Deltaxy[d1/3, slope1];** (*p4 is the point D*)
**slope2 = -1/slope1; p5 = p4 + Deltaxy[.25, slope2];** (*p5 is an arbitrary point on DH*)

(*To find the point H, use p4 and p5*): **x1 = p4[[1]]; y1 = p4[[2]]; x2 = p5[[1]]; y2 = p5[[2]];**
**dx = x2 - x1; dy = y2 - y1; dr = Sqrt[dx^2 + dy^2] ; D = x1*y2 - x2*y1; Del = Sqrt[dr^2 - D^2];**
**H1 = (D*dy + Sign[dy]*dx*Del)/dr^2; H2= (-D*dx + Abs[dy]*Del)/dr^2;**
**H = {H1,H2} ≈ { .5712240, .820794}** (*this is one of two solutions – at top and at bottom*)

(*Now find K and X*) **d2= EuclideanDistance[p1,H]; K = p1 - Deltaxy[d2, slope1];**
**d3 = EuclideanDistance[K, p4]; X = K + Deltaxy[d3/3, slope1] ≈ {0.4012520632830577, 0.5379964791784059}**

**m1 = (p1[[2]] - p3[[2]])/(p1[[1]] - p3[[1]]); m2 = (X[[2]] - p3[[2]])/(X[[1]] - p3[[1]]);**
**θ = N[ArcTan[(m2 - m1)/(1 + m1*m2)]] ≈ 0.5809057243292304 ≈ 33° 17' 0.40668"**
Compare to the actual which is **DMSString[100/3]** = 33° 20' 0"

For smaller angles, this method is so accurate that the secondary trisection can often be omitted and K can be used instead of X. Clearly K and H approach each other as the angle decreases.

**Example 2**: If Dürer had tested his method with the well-known angle of 60° it would have been impossible for him to discern any error as the following calculations show.

**p1 = {1,0}; p2= {Cos[60 Degree], Sin[60 Degree]};** H ≈ {0.9372458971632998, 0.3486690812942283}
X ≈ {0.8263539868031372, 0.3007637173887421};
Here m1 = 0 so **ArcTan[m2]** ≈ 0.34906100250062778410 = 19° 59' 59.00005"  Wow.

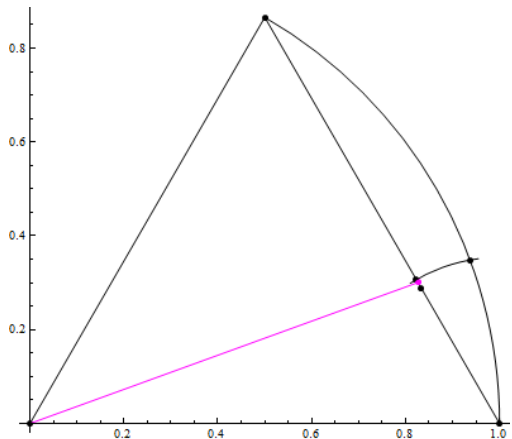

**Appendix B: Two 'modern' constructions**

If Dürer had gone one more step in his program for regular polygons he would have run into the 17-gon which was also known as the heptadecagon. Recall that the Greeks had no method for constructing regular prime polygons beyond N = 5 so their list of constructible regular polygons was 3, 4, 5, 6, 8, 10, 12, 15, 16, 20, 24, 30, 32, 40, 48, 60, 64, ...

These are all of the form: $2^m 3^k 5^j$ where m is a nonnegative integer and j and k are either 0 or 1.

That was where things stood for about 2000 years, when 19-year-old Carl Friedrich Gauss showed that a regular 17-gon was constructible in 1796. His solution to the cyclotomic equation $x^{17} = 1$ involved finding a sequence of nested quadratic equations based on the 8 pairs of matching vertices. This enabled him to obtain the following form for $\cos(2\pi/17)$

$$16\cos\frac{2\pi}{17} = -1+\sqrt{17}+\sqrt{34-2\sqrt{17}} +$$
$$2\sqrt{17+3\sqrt{17}-\sqrt{34-2\sqrt{17}}-2\sqrt{34+2\sqrt{17}}}.$$

Gauss showed 5 years later in *Disquisitiones Arithmeticae* that a prime regular polygon was constructible whenever it was a of the form $F_n = 2^{2^n} + 1$. These are called Fermat primes and the only known cases are $F_0 = 3$, $F_1 = 5$, $F_2 = 17$, $F_3 = 257$, $F_4 = 65537$. Gauss conjectured that the constructability condition was also necessary and this was proven by Pierre Wantzel in 1837. Therefore unless new Fermat primes are discovered (which is unlikely), there are only 5 prime regular polygons which can be constructed.

*A regular n-gon can be constructed with compass and straightedge iff n is the product of a power of 2 and any number of distinct Fermat primes.*

n = $2^m\, p_1 p_2..p_k$ where m is a non-negative integer and each $p_j$ is either 1 or the jth Fermat prime.

This yields the following (partial) list of constructible and non-constructible regular polygons:

| n | 3 | 4 | 5 | 6 | 7 | 8 | 9 | 10 | 11 | 12 | 13 | 14 | 15 | 16 | 17 |
|---|---|---|---|---|---|---|---|----|----|----|----|----|----|----|----|
| Constructible ? | Y | Y | Y | Y | N | Y | N | Y | N | Y | N | N | Y | Y | Y |

For more on the history of constructions, see the companion article Construction of Regular Polygons.pdf

Over the years there have been many different constructions for the regular pentagon and regular 17-gon. The constructions shown below are due to H.W. Richmond in 1893. (See *Mathematical Recreations and Essays* by W.W. Rouse Ball and H.M.S.Coxeter.) Both constructions are based on the same principle and they are simple and elegant. The pentagon construction takes just 2 steps and the 17-gon construction takes 6 steps.

For reference we have drawn the final version of the polygons as they would appear when inscribed in a unit circle, with vertex 1 at {1,0}. So in both diagrams P is at {1,0}.

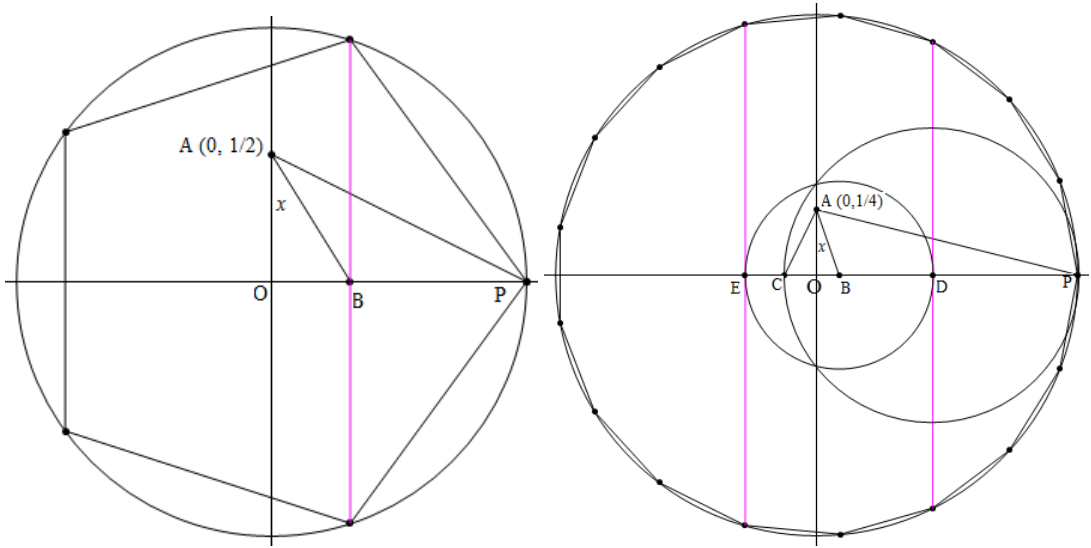

To construct the generating arc for the regular pentagon takes just 2 steps (on the left above)

    (i)      Locate point A at (0,1/2) as shown
    (ii)     Bisect angle OAP to obtain angle x and this defines point B and vertex 2.

To construct the generating arc for the regular 17-gon takes 6 steps

    (i)      Locate the point A at (0, ¼)
    (ii)     To find B, set the angle x to ¼ of angle OAP
    (iii)    To find C, make the angle CAB equal to Pi/4
    (iv)    The point D is the center of CP and this determines vertex 4 .
    (v)     Point E and vertex 6 are obtained from a circle centered at B passing through D
    (vi)    Now subtract the two arcs.

(If you are performing this 17-gon construction on a computer, the coordinates of B and C should be: B ≈ {{0.086037682852227701929, 0}, C ≈ {-0.12198209123162133118, 0})

In terms of the Lemoine complexity for constructions, there are less complex constructions using Carlyle circles, but no one knows the minimal complexity.

**On-Line Resources**

◆DürersHeptagon and Dürer9-gon are at demonstrations.wolfram.com. These are manipulates showing the construction of the 'regular' heptagon and 'regular' 9-gon. They are written in the new Mathematica computable data format.

◆There is an excellent data-base of Dürer's art at

http://www.wikipaintings.org/en/albrecht-durer/

◆The full text of the Four Books in German is available on-line at:

http://de.wikisource.org/wiki/Underweysung_der_Messung,_mit_dem_Zirckel_und_Richtscheyt,_in_Linien,_Ebenen_unnd_gantzen_corporen/Zweites_Buch

◆A copy of the 1538 edition is available at:

http://www.rarebookroom.org/Control/duruwm/index.html

◆A copy of the original 1525 edition is available at http://www.slub-dresden.de. This is the Saxon State Library.

◆Surprisingly, there is a copy of the 1525 edition available for sale at AbeBooks. The price is $52,250.00 (but shipping is just $10)

http://www.abebooks.com/servlet/SearchResults?an=Dürer&prl=10.00&recentlyadded=all&sortby=1&sts=t&x=55&y=9

◆Dave Richeson at Dickinson College has an excellent article on Durer's constructions at

http://divisbyzero.com/2011/03/22/albrecht-durers-ruler-and-compass-constructions/

◆DynamicsOfPolygons.org has two related articles which can be downloaded in pdf format:

(i) Tangent Durer.pdf which analyzes the dynamics of each non-regular 'Dürer polygon' under the outer billiards map (Tangent map)
(ii) Construction of Regular Polygons.pdf which discusses the general issue of compass and straightedge constructions using Gauss's *Disquisitiones Arithmeticae* of 1801.

## Bibliography

Ball W.W., Coxeter, H.M.S, *Mathematical Recreations and Essays*, Dover Publications, 1987

Conway, W.M. *Literary Remains of Albrecht Dürer – with transcripts from the British Museum Manuscripts,* Cambridge, 1889

Crowe, D.W. *Albrecht Dürer and the Regular Pentagon*, in Fivefold Symmetry by Istvan Hargittai, World Scientific, 1994

DeTemple, D.W. *Carlyle Circles and the Lemoine Simplicity of Polygon Constructions*, American Mathematical Monthly, Vol 98, No.2 (1991) pp. 97-108

Dürer, Albrecht, *Records of the Journey to Venice and the Low Countries*, edited by R. Fry, Boston, 1913

Dürer, Albrecht**,** *Underweysung der Messung mit dem Zirckel und Richtscheyt, Nuremberg,* 1525

Hunrath, K. *Albrecht Dürers Naherungskonstruktionen regelmassiger Vielecke*, Bibliotheca Mathematica, Dritte Folge, 1905

Pedoe, D., *Geometry and the Visual Arts*, Dover Publications, New York, 1976

Penofsky, E. *The Life and Art of Albrecht Dürer*, Princeton University Press, 1943

Richmonds, H.W. Quarterly Journal of Mathematics, 1893, vol 26, pg. 206

Staigmuller, H. *Dürer als Mathematiker, Programm des Koniglichen Realgymnasiums*, Stuttgart, 1891

Steck, Max*, Dürers Gestaltenlehre der Mathematik,* Halle, 1948

Strauss, W.L*. Albrecht Dürer -The Painters Manual,* Abaris Books, New York, 1977

The author of this paper is Emeritus Professor of Mathematics at California State University. He may be reached at ghhughes@csuchico.edu or mail@dynamicsofpolygons.org.